\documentclass[13pt]{article}
\usepackage{amsfonts}
\usepackage{amsmath}
\usepackage{amssymb}
\usepackage{amsthm}   
\usepackage{empheq}   
\usepackage{titlesec} 
\allowdisplaybreaks
\newtheorem{defi}{Definition}[section]
\newtheorem{thm}{Theorem}[section]
\newtheorem{lem}{Lemma}[section]
\newtheorem{prop}{Proposition}[section]
\newtheorem{rem}{Remark}[section]
\usepackage{bm}
\usepackage{bookmark}
\usepackage{arydshln}
\usepackage{lineno,hyperref}
\usepackage{xcolor}
\modulolinenumbers[5]

\begin{document}

\title{Stability and instability of the standing waves for the Klein-Gordon-Zakharov system in one space dimension}

\author{Silu Yin\thanks{yins11@shu.edu.cn, Department of Mathematics, Shanghai University, Shanghai 200444, P.R. China}}
\date{}

\maketitle

\begin{abstract}
  The orbital instability of standing waves for the Klein-Gordon-Zakharov system has been established in two and three space dimensions under radially symmetric condition, see Ohta-Todorova (SIAM J. Math. Anal. 2007). In the one space dimensional case, for the non-degenerate situation, we first check that the Klein-Gordon-Zakharov system satisfies Grillakis-Shatah-Strauss' assumptions on the stability and instability theorems for abstract Hamiltonian systems, see Grillakis-Shatah-Strauss (J. Funct. Anal. 1987). As to the degenerate case that the frequency $|\omega|=1/\sqrt{2}$, we follow Wu (ArXiv: 1705.04216, 2017) to describe the instability of the standing waves for the Klein-Gordon-Zakharov system, by using the modulation argument combining with the virial identity. For this purpose, we establish a modified virial identity to overcome several troublesome terms left in the traditional virial identity.
\end{abstract}
%


\linenumbers
\setcounter{section}{0}
\numberwithin{equation}{section}

\section{Introduction}
We study the Klein-Gordon-Zakharov (KGZ) system:
\begin{align}
&\partial_t^2 u-\partial_x^2 u+u+nu=0,&(t,x)\in\mathbb{R}\times\mathbb{R},\label{1.1}\\
&c_0^{-2}\partial_t^2 n-\partial_x^2 n=\partial_x^2(|u|^2),&(t,x)\in\mathbb{R}\times\mathbb{R},\label{1.2}
\end{align}
with the initial data
\begin{align}
  u(0,x)=u_0(x),\quad u_t(0,x)=v_0(x),\quad n(0,x)=n_0(x), \quad n_t(0,x)=\nu_0(x).
\end{align}
 It describes the nonlinear interaction between the quantum Langmuir and ion-acoustic waves in a plasma. Here $u$ is a complex valued function which represents the fast time scale component of an electric field raised by electrons, and $n$ is a real valued function which denotes the deviation of ion density from its equilibrium. For more physical backgrounds, the reader can refer to \cite{dendy,zakharov}.

The existence of the local smooth solutions to the Cauchy problem of the KGZ system can be proved by the standard Galerkin method (see, e.g., \cite{local}). Moreover, the global well-posedness was established in \cite{guo95,kato,kinoshita,m,ott}. Let $m:=-(-\Delta)^{-1}n_t$ with $m|_{t=0}:=-(-\Delta)^{-1}\nu_0$, where $\Delta=\partial_x^2$ in one space dimension. The solution $\vec{u}:=(u,u_t,n,n_t)^T=(u,v,n,\nu)^T$ of system \eqref{1.1}-\eqref{1.2} formally obeys the following energy, charge and momentum conservation identities:
\begin{align}
E(\vec{u})&=\frac12\|v\|^2+\frac1{4c_0^2}\| m_x\|^2+\frac12\|u_x\|^2+\frac12\|u\|^2+\frac14\|n\|^2+\frac12\int n|u|^2dx\notag\\
&=E(\vec{u}_0),\label{1.3}\\
Q(\vec{u})&=\textup{Im}\int \bar{u}vdx=Q(\vec{u}_0),\label{1.4}\\
P(\vec{u})&=2\textup{Re}\int u_t\bar{u}_xdx+\frac1{c_0^2}\int nm_xdx=P(\vec{u}_0),
\end{align}
 where $\vec{u}_0:=(u_0,v_0,n_0,\nu_0)^T$ and we denote
$$\|f\|:=\|f(x)\|_{L^2(\mathbb{R})}$$ for notational simplicity.
Without misunderstanding, we will use $X$ to represent the energy space $H^1(\mathbb{R})\times L^2(\mathbb{R})\times L^2(\mathbb{R})\times \dot{H}^{-1}(\mathbb{R})$ in this paper.

System \eqref{1.1}-\eqref{1.2} has the standing waves
$$(u_\omega(t,x),n_\omega(t,x))=(e^{i\omega t}\phi_\omega(x),-|\phi_\omega(x)|^2),$$
 where $\phi_\omega$ is the ground state solution to the following equation
\begin{equation}\label{1.5}
  -\partial_x^2\phi+(1-\omega^2)\phi-\phi^3=0.
\end{equation}
When the parameter $|\omega|<1$, equation \eqref{1.5} exists solution, see Strauss \cite{strauss}. In particular, the solution to \eqref{1.5} is unique.

In the present work, we are first interested in the stability and instability for the standing waves of the KGZ system with non-degenerate frequency $\omega\in(-1,1)$.
We introduce the following notation for $\vec{f}=(f,g,h,k)^T\in X$:
$$T(\theta)\vec{f}=(e^{i\theta}f,e^{i\theta}g,h,k)^T.$$
Let
$$\vec{\Phi}_\omega=(\phi_\omega,i\omega\phi_\omega,-|\phi_\omega|^2,0)^T,$$
and denote
$$U_\varepsilon(\vec{\Phi}_\omega)=\{\vec{u}\in X:\inf_{(\theta,y)\in\mathbb{R}\times\mathbb{R}}\|\vec{u}-T(\theta)\vec{\Phi}_\omega(\cdot-y)\|_X<\varepsilon\}.$$

\begin{defi}
The standing wave solution $(u_\omega,n_\omega)$ of the KGZ system \eqref{1.1}-\eqref{1.2} is said to be orbitally stable if for any given $\varepsilon>0$, there exists $\delta>0$, such that when $\vec{u}_0\in U_\delta(\vec{\Phi}_\omega)$, the solution $\vec{u}(t)$ of system \eqref{1.1}-\eqref{1.2} with the initial data $\vec{u}_0$ exists for all $t>0$, and $\vec{u}(t)\in U_\varepsilon(\vec{\Phi}_\omega)$ for all $t>0$. Otherwise, $(u_\omega,n_\omega)$ is said to be orbitally unstable.
\end{defi}

 It is worthy to mention that Grillakis-Shatah-Strauss \cite{grillakis1,grillakis2} established the sability and instability theory for abstract Hamiltonian systems. Consider the KGZ system \eqref{1.1}-\eqref{1.2} in the Hamiltonian case
\begin{equation}
\frac{d\vec{u}}{dt}=JE'(\vec{u}),
\end{equation}
where $E'$ is the Fr\'{e}chet derivative of $E$, and define the functional $$S_\omega(\vec{u})=E(\vec{u})-\omega Q(\vec{u}).$$
We can verify that $\vec{\Phi}_\omega$ is a critical point of $S_\omega(\vec{u})$. Via analyzing the spectrum, we can obtain the orbital stability according to the framework of Grillakis-Shatah-Strauss \cite{grillakis1} in the non-degenerate case. In other words, if $S''_\omega(\vec{u})$ has exactly one negative eigenvalue, then it is orbitally stable if and only if $S_\omega(\vec{u})$ is strictly convex at $\omega$. Specifically, by checking the assumptions of the stability and instability theorems in \cite{grillakis1}, we can get directly the following consequence.
\begin{prop}\label{thm0}
The standing wave $(e^{i\omega t}\phi_\omega,-|\phi_\omega|^2)$ of the KGZ system \eqref{1.1}-\eqref{1.2} is orbitally unstable for $|\omega|<\frac1{\sqrt{2}}$ and stable for $\frac1{\sqrt{2}}<|\omega|<1$.
\end{prop}


For the case $|\omega|=\frac1{\sqrt{2}}$, which is degenerate based on \cite{grillakis1,grillakis2}, Comech-Pelinovsky \cite{comech} showed the instability by a careful analysis of the linearized system. Ohta gave a different and shorter proof in \cite{ohta} by the Lyapunov functional method. Maeda \cite{maeda} extended the previously mentioned researches to subtler assumptions of $S_\omega$ based on a purely variational argument. Here, we will use the modulation argument combining with the virial identity to obtain the instability result. This method is inspired by the work of Wu \cite{w}, in which he established the orbital instability of the nonlinear Klein-Gordon equation in the critical frequency case in one space dimension. The modulation argument was introduced by Weinstein \cite{weinstein} for nonlinear Schr\"{o}dinger equations, developed and widely applied by Martel-Merle \cite{martel},  Martel-Merle-Tsai \cite{mm}, Merle \cite{merle}, Merle-Rapha\"{e}l \cite{raphael}, Stuart \cite{stuart}, Bellazzini-Ghimenti-Le Coz \cite{b}, Le Coz-Wu \cite{c}. The local version of the virial type identity is crucial in the proof of our main result. For the KGZ system, there will be several troublesome terms left in the original virial identity, so we need a modified virial identity, see Section 3. Our main result is
\begin{thm}\label{thm}
When $|\omega|=\frac1{\sqrt{2}}$, the standing wave $(e^{i\omega t}\phi_\omega,-|\phi_\omega|^2)$ of the KGZ system \eqref{1.1}-\eqref{1.2} is orbitally unstable.
\end{thm}

Before ending this introduction, let us mention some related work in several space dimensions. Ohta-Todorova \cite{ot} showed that the radially symmetric standing wave $ (e^{i\omega t}\phi_\omega(x),-|\phi_\omega(x)|^2)$ of the KGZ system is strongly unstable in two and three space dimensions when $\omega\in(-1,1)$. One key step in their proof is that the decay estimate for radially symmetric functions in $H^1(\mathbb{R}^d)$,
$$\|u(x)\|\leq C|x|^{(1-d)/2}\|u\|_{H^1(\mathbb{R}^d)}, \quad\textup{for}\ |x|\geq1$$
is valid when $d\geq2$ but failed when $d=1$, see \cite{strauss}. The case of standing wave with the ground state ($\omega=0$) was discussed by Gan-Zhang \cite{gan} in three space dimensions. For more related works in two and higher space dimensions, one can see \cite{gan09,glangetas,guozihua,hakkaev} and the references therein.

The plan of this paper is as follows. In Section 2, we give some preliminaries, in which we present several basic properties and the coercivity of the Hessian. In Section 3, we introduce the modified virial type identity. In Section 4, we explain the modulation theory in a neighborhood of the standing wave, and control the terms of the modified virial identity. Finally, we prove the main theorem in Section 5.

\section{Basic properties}
We first infer that from the definition of $E(\vec{u})$, we have
\begin{equation}
E'(\vec{u})=\begin{pmatrix}
                                         -\partial_x^2 u+u+nu \\
                                         v \\
                                         \frac12n+\frac12|u|^2 \\
                                         \frac1{2c_0^2}(-\Delta)^{-1}\nu \\
                                       \end{pmatrix}.
\end{equation}
Let
$$J=\begin{pmatrix}
                            0  & 1 & 0 & 0 \\
                            -1 & 0 & 0 & 0 \\
                            0 & 0& 0 & -2c_0^2\partial_x^2 \\
                            0 & 0 & 2c_0^2\partial_x^2 & 0 \\
                          \end{pmatrix},$$
which is a skew symmetric linear operator. We can rewrite the KGZ system in the form of
\begin{equation}\label{3.16}
\frac{d\vec{u}}{dt}=JE'(\vec{u}).
\end{equation}
Define the inner product
$$\langle \vec{f},\vec{g}\rangle=\textup{Re}\int_\mathbb{R}\vec{f}(x)^T\cdot\overline{\vec{g}(x)} dx.$$
In what follows, we write $A\gtrsim B$ (or $A\lesssim B$) to express $A\geq kB$ (or $A\leq k B$) for certain positive constant $k>0$. We have the following basic properties.

\begin{lem}\label{lem3.1}
For the KGZ system \eqref{1.1}-\eqref{1.2}, we have
  \begin{equation}\label{3.6-0}
    \frac d{dw}Q(\vec{\Phi}_\omega)=\frac{1-2\omega^2}{\sqrt{1-\omega^2}}\|\phi_0\|^2
  \end{equation}
  and
  \begin{equation}\label{3.7}
    3E(\vec{\Phi}_\omega)-4\omega Q(\vec{\Phi}_\omega)=(1-2\omega^2)\|\phi_\omega\|^2=(1-2\omega^2)\sqrt{1-\omega^2}\|\phi_0\|^2.
  \end{equation}
\end{lem}
\textit{Proof.} From the definition of the charge and rescaling, we have
$$\phi_\omega(x)=\sqrt{1-\omega^2}\phi_0(\sqrt{1-\omega^2}x),$$
in which $\phi_0(x)$ is the solution of
\begin{equation}
  -\partial_x^2\phi+\phi-\phi^3=0.
\end{equation}
Then we get
\begin{align}\label{2.3}
  Q(\vec{\Phi}_\omega)=\omega\|\phi_\omega\|^2=\omega\sqrt{1-\omega^2}\|\phi_0\|^2,
\end{align}
then \eqref{3.6-0} can be obtained by taking the derivative of \eqref{2.3} with respect to $\omega$.

Let $\phi_\omega$ and $x\partial_x\phi_\omega$ be the multipliers of \eqref{1.5}, respectively. We have
\begin{align*}
  \|\partial_x\phi_\omega\|^2+(1-\omega^2)\|\phi_\omega\|^2-\|\phi_\omega\|_{L^4}^4=0
\end{align*}
and
\begin{align*}
  \|\partial_x\phi_\omega\|^2-(1-\omega^2)\|\phi_\omega\|^2+\frac12\|\phi_\omega\|_{L^4}^4=0.
\end{align*}
Solving the above equations, we get
$$\|\partial_x\phi_\omega\|^2=\frac{1-\omega^2}3\|\phi_\omega\|^2,\quad \|\phi_\omega\|_{L^4}^4=\frac{4(1-\omega^2)}3\|\phi_\omega\|^2,$$
so
\begin{align*}
  3E(\vec{\Phi}_\omega)-4\omega Q(\vec{\Phi}_\omega)=&\frac32\|\partial_x\phi_\omega\|^2+\frac{3(1+\omega^2)}2\|\phi_\omega\|^2-\frac34\|\phi_\omega\|_{L^4}^4-4\omega^2\|\phi_\omega\|^2\\
  =&(1-2\omega^2)\|\phi_\omega\|^2=(1-2\omega^2)(1-\omega^2)^\frac12\|\phi_0\|^2.
\end{align*}
This proves \eqref{3.7}.
\hfill$\Box$

Define the functional
\begin{align}\label{27}
  S_\omega(\vec{u})=E(\vec{u})-\omega Q(\vec{u}).
\end{align}
We write its derivative as $\langle S'_\omega(\vec{u}),\vec{g}\rangle$, where $S'_\omega$ is a functional from $X$ to its dual, and its second derivative as $\langle S''_\omega(\vec{u})\vec{f},\vec{g}\rangle$, for $\vec{g}\in X$.

Let \begin{align}
  \vec{\Upsilon}_\omega&=(i\phi_\omega,-\omega\phi_\omega,0,0)^T,\\ \partial_x\vec{\Phi}_\omega&=(\partial_x\phi_\omega,i\omega\partial_x\phi_\omega,-2\phi_\omega\partial_x\phi_\omega,0)^T,\\
  \vec{\Psi}_\omega&=(2\omega\phi_\omega,0,0,0)^T.
\end{align}
We will give the coercivity property in Lemma \ref{lem2.3} below. Before reaching this, we show some preliminaries.
\begin{lem}\label{lem2.2}
For the KGZ system \eqref{1.1}-\eqref{1.2}, we have
\begin{equation}\label{2.4-0}
 S'_\omega(\vec{\Phi}_\omega)=0
\end{equation}
and
\begin{equation}\label{2.5-00}
  \textup{Ker}(S''_\omega(\vec{\Phi}_\omega))=\textup{Span}\{\vec{\Upsilon}_\omega,\partial_x\vec{\Phi}_\omega\},
\end{equation}
where $S_\omega$ is defined by \eqref{27}.
\end{lem}
\textit{Proof.}
By definitions of the energy and charge, we have
\begin{equation}
  S'_\omega(\vec{u})=E'(\vec{u})-\omega Q'(\vec{u})=\begin{pmatrix}
                       -u_{xx}+u+nu \\
                       v \\
                       \frac12n+\frac12|u|^2 \\
                       \frac1{2c_0^2}(-\Delta)^{-1}\nu \\
                     \end{pmatrix}+i\omega\begin{pmatrix}
                                            v \\
                                            -u\\
                                            0 \\
                                            0 \\
                                          \end{pmatrix}
\end{equation}
and
\begin{equation}\label{2.5-0}
  S''_\omega(\vec{\Phi}_\omega)\vec{f}=\begin{pmatrix}
                       -f_{xx}+f-\phi_\omega^2f+\phi_\omega h \\
                       g \\
                       \frac12h+\phi_\omega\textup{Re}f \\
                       \frac1{2c_0^2}(-\Delta)^{-1}k \\
                     \end{pmatrix}+i\omega\begin{pmatrix}
                                            g \\
                                            -f\\
                                            0 \\
                                            0 \\
                                          \end{pmatrix}.
\end{equation}
Then we have
\begin{equation*}
 S'_\omega(\vec{\Phi}_\omega)=E'(\vec{\Phi}_\omega)-\omega Q'(\vec{\Phi}_\omega)=0
\end{equation*}
and
\begin{align}\label{2.11}
S''_\omega(\vec{\Phi}_\omega)\vec{\Upsilon}_\omega=S''_\omega(\vec{\Phi}_\omega)\partial_x\vec{\Phi}_\omega=0.
\end{align}
It follows from \eqref{2.11} that
\begin{align}\label{16}
\textup{Span}\{\vec{\Upsilon}_\omega,\partial_x\vec{\Phi}_\omega\}\subset  \textup{Ker}(S''_\omega(\vec{\Phi}_\omega)).
\end{align}

Next we show that for any given $\vec{f}=(f,g,h,k)\in  \textup{Ker}(S''_\omega(\vec{\Phi}_\omega))$ with $\vec{f}\in X$, we have $\vec{f}\in\textup{Span}\{\vec{\Upsilon}_\omega,\partial_x\vec{\Phi}_\omega\}$.
In fact, from \eqref{2.5-0}, if $\vec{f}=(f,g,h,k)\in  \textup{Ker}(S''_\omega(\vec{\Phi}_\omega))$ and $\vec{f}\neq 0$, we have
\begin{align}
\left\{
\begin{array}{rrrr}
 -f_{xx}+f-\phi_\omega^2f+\phi_\omega h+i\omega g&=0, \\
                       g-i\omega f&=0 ,\\
                       \frac12h+\phi_\omega\textup{Re}f&=0, \\
                       \frac1{2c_0^2}(-\Delta)^{-1}k&=0.
\end{array}\right.
\end{align}
Thus, $\vec{f}=(f,i\omega f,-2\phi_\omega\textup{Re}f,0)$, and $f$ solves the equation
\begin{align}\label{2.9-0}
  -f_{xx}+(1-\omega^2)f-\phi_\omega^2f-2\phi_\omega^2 \textup{Re}f=0.
\end{align}
Let
$$f=f_1+i f_2,\quad f_1=\textup{Re}f,\quad f_2=\textup{Im}f$$
and
\begin{align}
  L_+&=-\partial_x^2+(1-\omega^2)-3\phi_\omega^2,\\
  L_-&=-\partial_x^2+(1-\omega^2)-\phi_\omega^2.
\end{align}
Equation \eqref{2.9-0} implies that
$$f_1\in\textup{Ker}(L_+),\quad f_2\in\textup{Ker}(L_-).$$
In order to show that $f_1$ has a linear correlation with $\partial_x\phi_\omega$, we consider their Wronskian determinant:
\begin{align}
  W_+=\begin{vmatrix}
        f_1 & \partial_x\phi_\omega \\
        \partial_x f_1 & \partial_x^2\phi_\omega \\
      \end{vmatrix}=f_1\partial_x^2\phi_\omega-\partial_x f_1\partial_x\phi_\omega.
\end{align}
Since $\partial_x\phi_\omega,f_1\in\textup{Ker}(L_+)$, we have
\begin{align*}
  \partial_x^3\phi_\omega&=(1-\omega^2)\partial_x\phi_\omega-3\phi_\omega^2\partial_x\phi_\omega,\\
  \partial_x^2f_1&=(1-\omega^2)f_1-3\phi_\omega^2f_1.
\end{align*}
Therefore
\begin{equation}\label{22}
  \begin{split}
    W'_+&=f_1\partial_x^3\phi_\omega-\partial_x^2 f_1\partial_x\phi_\omega\\
    &=f_1(1-\omega^2)\partial_x\phi_\omega-3f_1\phi_\omega^2\partial_x\phi_\omega-(1-\omega^2)f_1\partial_x\phi_\omega+3\phi_\omega^2 f_1\partial_x\phi_\omega\\
    &=0.
  \end{split}
\end{equation}
By the exponential decay of $\phi_\omega$, which is well-known (see \cite{gidas,kwong}), we have
\begin{align}
  \lim_{x\rightarrow\pm\infty}W_+=\lim_{x\rightarrow\pm\infty}(f_1\partial_x^2\phi_\omega-\partial_x f_1\partial_x\phi_\omega)=0.
\end{align}
So we get that $W_+=0$ with \eqref{22}.
Thus we have
\begin{align}\label{2.15}
  f_1\in\textup{Span}\{\partial_x\phi_\omega\}.
\end{align}
For $L_-$, we analogously consider the Wronskian determinant of $f_2$ and $\phi_\omega$:
\begin{align}
  W_-=\begin{vmatrix}
        f_2 & \phi_\omega \\
        \partial_x f_2 & \partial_x\phi_\omega \\
      \end{vmatrix}=f_2\partial_x\phi_\omega-\partial_x f_2\phi_\omega,
\end{align}
and then
\begin{equation}
  \begin{split}
    W'_-&=f_2\partial_x^2\phi_\omega-\partial_x^2f_2\phi_\omega\\
    &=f_2[(1-\omega^2)\phi_\omega-\phi_\omega^3]-[(1-\omega^2)f_2-\phi_\omega^2f_2]\phi_\omega\\
    &=0.
  \end{split}
\end{equation}
So we have
$W_-=0$, which implies
\begin{align}\label{2.18}
  f_2\in\textup{Span}\{\phi_\omega\}.
\end{align}
It follows from \eqref{2.15} and \eqref{2.18} that $ f\in\textup{Span}\{\partial_x\phi_\omega,i\phi_\omega\}$. This implies that
\begin{align}
\textup{Ker}(S''_\omega(\vec{\Phi}_\omega))\subset\textup{Span}\{\vec{\Upsilon}_\omega,\partial_x\vec{\Phi}_\omega\},
\end{align}
which combined with \eqref{16} gives \eqref{2.5-00}.
\hfill$\Box$

The following lemma tells us that $S''_\omega(\vec{\Phi}_\omega)$ has a negative direction.
\begin{lem}\label{lem2.3-0}
  Let $\vec{F}_\omega=(\partial_\omega\phi_\omega,i\omega\partial_\omega\phi_\omega,-2\phi_\omega\partial_\omega\phi_\omega,0)^T$. We have
  \begin{equation}\label{2.4}
   S''_\omega(\vec{\Phi}_\omega)\vec{F}_\omega=\vec{\Psi}_\omega.
  \end{equation}
Then, if $\omega\neq 0$, we have
  \begin{equation}\label{2.5}
    \langle S''_\omega(\vec{\Phi}_\omega)\vec{F}_\omega,\vec{F}_\omega\rangle <0.
  \end{equation}
\end{lem}
\textit{Proof.}
Since
$$ -\partial_x^2\phi_\omega+(1-\omega^2)\phi_\omega-\phi_\omega^3=0,$$
we get
\begin{equation}
  -\partial_x^2\partial_\omega\phi_\omega+(1-\omega^2)\partial_\omega\phi_\omega-3\phi_\omega^2\partial_\omega\phi_\omega=2\omega\phi_\omega.
\end{equation}
Thus
\begin{equation}
  S''_\omega(\vec{\Phi}_\omega)\vec{F}_\omega=\begin{pmatrix}
                      -\partial_x^2\partial_\omega\phi_\omega+\partial_\omega\phi_\omega-\phi_\omega^2\partial_\omega\phi_\omega-2\phi_\omega^2\partial_\omega\phi_\omega-\omega^2\partial_\omega\phi_\omega \\
                      i\omega\partial_\omega\phi_\omega -i\omega\partial_\omega\phi_\omega  \\
                       -\phi_\omega\partial_\omega\phi_\omega+\phi_\omega\partial_\omega\phi_\omega \\
                       0 \\
                     \end{pmatrix}=\begin{pmatrix}
                                           2\omega\phi_\omega \\
                                            0\\
                                            0 \\
                                            0 \\
                                          \end{pmatrix}.
\end{equation}
By rescaling, we have
\begin{align}
  \langle S''_\omega(\vec{\Phi}_\omega)\vec{F}_\omega,\vec{F}_\omega\rangle=\omega\partial_\omega(\|\phi_\omega\|^2)=\omega\partial_\omega[(1-\omega^2)^\frac12\|\phi_0\|^2]=-\frac{\omega^2}{\sqrt{1-\omega^2}}\|\phi_0\|^2<0
\end{align}
for $\omega\neq 0$. \hfill$\Box$

Now we are ready to present the following coercivity property. Similar argument can be found in \cite{b} for Klein-Gordon equation and in \cite{c} for nonlinear Schr\"{o}dinger equation.
\begin{lem}[Coercivity]\label{lem2.3}
  Assume $\omega\neq 0$, for any given $\vec{\xi}=(\xi,\eta,\zeta,\iota)^T\in X$ satisfying
  \begin{align}\label{2.26}
    \langle \vec{\xi},\vec{\Upsilon}_\omega\rangle=\langle\vec{\xi},\partial_x\vec{\Phi}_\omega\rangle=\langle\vec{\xi},\vec{\Psi}_\omega\rangle=0,
  \end{align}
  there exists a positive constant $\delta$ such that
\begin{align}\label{35}
  \langle S''_\omega(\vec{\Phi}_\omega)\vec{\xi},\vec{\xi}\rangle\geq \delta \|\vec{\xi}\|_X^2.
\end{align}
\end{lem}
\textit{Proof.}
\textit{Step 1: Spectral analysis.}

Let
\begin{equation}
\mathcal{L}:=\begin{pmatrix}
  -\partial_x^2+1 & i\omega& 0& 0 \\
  -i\omega & 1 & 0 & 0 \\
  0 & 0 & \frac12 & 0 \\
  0 & 0 & 0& \frac1{2c_0^2}(-\Delta)^{-1} \\
\end{pmatrix}.
\end{equation}
For $\vec{f}=(f,g,h,k)\in X$, we have
\begin{align}
  \langle \mathcal{L}f,f\rangle=\|f_x \|^2+\|f\|^2+\|g\|^2+\frac12\|h\|^2+\frac1{2c_0^2}\|(-\Delta)^{-1}k_x\|^2\geq \delta\|\vec{f}\|_{X}^2,
\end{align}
where
\begin{align}\label{2.29-0}
  \delta=\min\{\frac12,\frac1{2c_0^2}\}>0.
\end{align}
This means that the essential spectrum of $\mathcal{L}$ is positive and away from zero. Since $S''_\omega(\vec{\Phi}_\omega)$ is a compact perturbation of $\mathcal{L}$, the essential spectrum of $S''_\omega(\vec{\Phi}_\omega)$ is also positive and away from zero by Weyl's Theorem. The rest of its spectrum consists of isolated eigenvalues. By the variational characterization of $\vec{\Phi}_\omega$, see e.g. \cite{ambrosetti}, the Morse Index of $S_\omega$ is $1$ at most. Then combined with Lemma \ref{lem2.3-0}, we infer that $S''_\omega(\vec{\Phi}_\omega)$ admits only one negative eigenvalue $-\gamma<0$ satisfying
\begin{align}\label{2.29}
  S''_\omega(\vec{\Phi}_\omega)\Gamma=-\gamma\Gamma
\end{align}
with $\|\Gamma\|_X=1$.

\textit{Step 2: Positivity analysis.}

We claim that if $\vec{\xi}$ satisfies  the orthogonality condition \eqref{2.26}, then
\begin{align}
 \langle S''_\omega(\vec{\Phi}_\omega)\vec{\xi},\vec{\xi}\rangle>0.
\end{align}
In fact, by Lemma \ref{lem2.2} and \eqref{2.26}, we can write the orthogonal decomposition of $\vec{\xi}$ along the spectrum of $S''_\omega(\vec{\Phi}_\omega)$ as follows:
\begin{align}
  \vec{\xi}=\alpha\vec{\Gamma}+\vec{\Theta},
\end{align}
where $\alpha$ is a real number and $\vec{\Theta}$ is in the positive eigenspace of $S''_\omega(\vec{\Phi}_\omega)$. If $\alpha=0$, the desired conclusion is trivial.

 For $\alpha\neq 0$, 
We decompose $\vec{F}_\omega$ orthogonally  along the spectrum of $S''_\omega(\vec{\Phi}_\omega)$:
\begin{align}
 \vec{F}_\omega=\beta\vec{\Gamma}+\vec{\Omega}+\vec{\Pi},
\end{align}
where $\beta\neq0$, $\vec{\Omega}$ is in the positive eigenspace of $S''_\omega(\vec{\Phi}_\omega)$ and $\vec{\Pi}\in \textup{Ker}(S''_\omega(\vec{\Phi}_\omega))$. In particular, we have
\begin{align}
 \langle S''_\omega(\vec{\Phi}_\omega)\Omega,\Omega\rangle\gtrsim \|\Omega\|_X.
\end{align}
Because of \eqref{2.26}, we have
\begin{equation}
  \begin{split}
    0&=  \langle \vec{\Psi}_\omega,\vec{\xi}\rangle\\
    &=\langle S''_\omega(\vec{\Phi}_\omega)\vec{F}_\omega,\vec{\xi}\rangle\\
    &=\langle S''_\omega(\vec{\Phi}_\omega)(\beta\vec{\Gamma}+\vec{\Omega}+\vec{\Pi}),\alpha\vec{\Gamma}+\vec{\Theta}\rangle\\
    &=-\gamma \alpha\beta+\langle S''_\omega(\vec{\Phi}_\omega)\vec{\Omega},\vec{\Theta}\rangle.
  \end{split}
\end{equation}
Noting that both $\vec{\Omega}$ and $\vec{\Theta}$ are in the positive eigenspace of $S''_\omega(\vec{\Phi}_\omega)$, by Cauchy-Schwarz inequality and the fact that $\langle S''_\omega(\vec{\Phi}_\omega)\vec{\Omega},\vec{\Theta}\rangle=\langle S''_\omega(\vec{\Phi}_\omega)\vec{\Theta},\vec{\Omega}\rangle$, we have
\begin{align}
  \gamma^2 \alpha^2\beta^2=\langle S''_\omega(\vec{\Phi}_\omega)\vec{\Omega},\vec{\Theta}\rangle^2\leq\langle S''_\omega(\vec{\Phi}_\omega)\vec{\Omega},\vec{\Omega}\rangle\langle S''_\omega(\vec{\Phi}_\omega)\vec{\Theta},\vec{\Theta}\rangle.
\end{align}
Moreover, \eqref{2.5} tells us that
\begin{equation}
\begin{split}
  \langle S''_\omega(\vec{\Phi}_\omega)\vec{F}_\omega,\vec{F}_\omega\rangle&=\langle S''_\omega(\vec{\Phi}_\omega)(\beta\vec{\Gamma}+\vec{\Omega}+\vec{\Pi}),\beta\vec{\Gamma}+\vec{\Omega}+\vec{\Pi}\rangle\\
  &=-\gamma \beta^2+\langle S''_\omega(\vec{\Phi}_\omega)\vec{\Omega},\vec{\Omega}\rangle<0,
  \end{split}
\end{equation}
Therefore,
\begin{equation}
  \begin{split}
     \langle S''_\omega(\vec{\Phi}_\omega)\vec{\xi},\vec{\xi}\rangle&=\langle S''_\omega(\vec{\Phi}_\omega)(\alpha\vec{\Gamma}+\vec{\Theta}),\alpha\vec{\Gamma}+\vec{\Theta}\rangle\\
     &=-\gamma \alpha^2+\langle S''_\omega(\vec{\Phi}_\omega)\vec{\Theta},\vec{\Theta}\rangle\\
     &\geq -\gamma \alpha^2+\frac{\gamma^2 \alpha^2 \beta^2}{\langle S''_\omega(\vec{\Phi}_\omega)\vec{\Omega},\vec{\Omega}\rangle}\\
     &> -\gamma \alpha^2+\frac{\gamma^2 \alpha^2 \beta^2}{\gamma \beta^2}=0.
  \end{split}
\end{equation}

\textit{Step 3: Coercivity analysis.} We prove the coercivity by contradiction. Assume that $\vec{\xi}_n\ (n\in\mathbb{N}^+)$ satisfy
  \begin{align}
    \langle \vec{\xi}_n,\vec{\Upsilon}_\omega\rangle=\langle\vec{\xi}_n,\partial_x\vec{\Phi}_\omega\rangle=\langle\vec{\xi}_n,\vec{\Psi}_\omega\rangle=0
  \end{align}
 with $\|\vec{\xi}_n\|_X=1$ and
 \begin{align}\label{2.39}
   \lim_{n\rightarrow\infty}\langle S''_\omega(\vec{\Phi}_\omega)\vec{\xi}_n,\vec{\xi}_n\rangle=0,
 \end{align}
then there exists $\vec{\xi}_\infty\in X$, such that $\{\vec{\xi}_n\}$ converges to $\vec{\xi}_\infty$ weakly and $\vec{\xi}_\infty$ verifies
   \begin{align}
    \langle \vec{\xi}_\infty,\vec{\Upsilon}_\omega\rangle=\langle\vec{\xi}_\infty,\partial_x\vec{\Phi}_\omega\rangle=\langle\vec{\xi}_\infty,\vec{\Psi}_\omega\rangle=0.
  \end{align}

 We first prove that
 \begin{equation}\label{50}
   \vec{\xi}_\infty=0.
 \end{equation}
 Otherwise, from the positivity analysis in Step 2, we have
   \begin{align}
  \langle S''_\omega(\vec{\Phi}_\omega)\vec{\xi}_\infty,\vec{\xi}_\infty\rangle>0.
 \end{align}
Moreover, by the weak convergence and exponential decay of $\vec{\Phi}_\omega$, we infer that
  \begin{align}
    \langle S''_\omega(\vec{\Phi}_\omega)\vec{\xi}_\infty,\vec{\xi}_\infty\rangle\leq\lim_{n\rightarrow\infty}\langle S''_\omega(\vec{\Phi}_\omega)\vec{\xi}_n,\vec{\xi}_n\rangle=0.
  \end{align}
  So we have \eqref{50}.

  On the other hand, it follows from \eqref{2.5-0} and $\|\vec{\xi}_n\|_X=1$ that
  \begin{equation}\label{2.43}
  \begin{split}
     \langle S''_\omega(\vec{\Phi}_\omega)\vec{\xi}_n,\vec{\xi}_n\rangle=&\|\partial_x\xi_n\|^2+\|\xi_n\|^2+\|\eta_n\|^2+\frac12\|\zeta_n\|^2+\frac1{2c_0^2}\|(-\Delta)^{-1}\partial_x\iota\|^2\\
     &-\|\phi_\omega\xi_n\|^2+2\langle\phi_\omega\zeta_n,\xi_n\rangle-2\langle i\omega\xi_n,\zeta_n\rangle\\
     \geq&\delta-\|\phi_\omega\xi_n\|^2+2\langle\phi_\omega\zeta_n,\xi_n\rangle-2\langle i\omega\xi_n,\zeta_n\rangle,
  \end{split}
\end{equation}
where $\delta$ is defined by \eqref{2.29-0}.
Taking $n\rightarrow\infty$ in \eqref{2.43}, it follows from \eqref{2.39} that
  \begin{equation}
  \begin{split}
   -\|\phi_\omega\xi_\infty\|^2+2\langle\phi_\omega\zeta_\infty,\xi_\infty\rangle-2\langle i\omega\xi_\infty,\zeta_\infty\rangle\leq-\delta<0.
  \end{split}
\end{equation}
This gives a contradiction with  \eqref{50}. Hence, for any given $\vec{\xi}$ satisfying \eqref{2.26}, $\langle S''_\omega(\vec{\Phi}_\omega)\vec{\xi},\vec{\xi}\rangle$ is not only positive, but also strictly bigger than a positive constant. This proves the coercivity property.  \hfill$\Box$

\begin{rem}
When $\omega\neq0$, Lemma \ref{lem2.3-0} and Lemma \ref{lem2.3} imply that $S''_\omega$ has exactly one negative eigenvalue and its essential spectrum is positive and bounded away from zero. From Lemma \ref{lem3.1}-Lemma \ref{lem2.3}, we infer that the KGZ system obeys Grillakis-Shatah-Strauss's assumptions for abstract Hamiltonian systems. When $\omega=0$, $(\phi_0,-|\phi_0|^2)$ is a stationary solution of the KGZ system. Grillakis-Shatah-Strauss gave the sharp conditions for the stability and instability of solitary waves in \cite{grillakis1} as follows:
 \begin{itemize}
   \item If $-\frac{d}{d\omega}Q(\phi_\omega)<0$, then the  $\phi$-orbit is unstable;
   \item The fact that $\phi$-orbit is stable if and only if $-\frac{d}{d\omega}Q(\phi_\omega)>0$.
 \end{itemize}
Hence, the standing wave of the KGZ system is unstable for $|\omega|<\frac1{\sqrt{2}}$ and stable for $\frac1{\sqrt{2}}<|\omega|<1$.
\end{rem}


\section{Modified Virial Type Identity}
In the degenerate case $|\omega|=\frac1{\sqrt{2}}$, the proof of Theorem \ref{thm} is based on the local versions of the virial type identity. We first define a cutoff function $\varphi_R\in C^\infty(\mathbb{R})$ by
\begin{equation}
  \varphi_R(x)=\left\{\begin{array}{ll}
  x,\quad |x|\leq R,\\
    0,\quad |x|\geq 2R.
  \end{array}\right.
\end{equation}
\begin{lem}\label{lem3-1}
Let the traditional virial type quantity $I(t)$ be defined by
\begin{align}\label{3.0}
  I(t)=\textup{Re}\int_\mathbb{R} u\bar{u}_tdx+2\textup{Re}\int_\mathbb{R}\varphi_R(x-y(t)) u_x\bar{u}_tdx+\frac1{c_0^2}\int_\mathbb{R}\varphi_R(x-y(t))nm_xdx.
\end{align}
 Then for the KGZ system, we have
\begin{equation}\label{3.5}
\begin{split}
I'(t)=&-\dot{y}P(\vec{u}_0)+\dot{y}\int_\mathbb{R}[1-\varphi'_R(x-y(t))][2\textup{Re}(u_x\bar{u}_t)+\frac1{c_0^2}nm_x]dx\\
&-2\|u_x\|^2-\frac12\|n\|^2-\frac1{2c_0^2}\|m_x\|^2-\int_\mathbb{R}n|u|^2dx\\
&+\int_\mathbb{R}[1-\varphi'_R(x-y(t))][|u_t|^2+|u_x|^2-|u|^2+\frac1{2c_0^2}|m_x|^2+\frac12|n|^2]dx.
\end{split}
\end{equation}
\end{lem}
{\it Proof.} By elementary calculation, we have
\begin{align}
\frac d{dt}\textup{Re}\int_\mathbb{R} u\bar{u}_tdx&=\|u_t\|^2-\|u_x\|^2-\|u\|^2-\int_\mathbb{R} n|u|^2dx,\label{3.2}\\
2\textup{Re}\int_\mathbb{R}\varphi_R\frac d{dt}(u_x\bar{u}_t)dx&=-\int_\mathbb{R}\varphi'_R|u_t|^2dx-\int_\mathbb{R}\varphi'_R|u_x|^2dx+\int_\mathbb{R}\varphi'_R|u|^2dx+\int_\mathbb{R}(\varphi_Rn)_x|u|^2dx
\end{align}
and
\begin{align}
\frac1{c_0^2}\int_\mathbb{R}\varphi_R\frac d{dt}(nm_x)dx=-\frac1{2c_0^2}\int_\mathbb{R}\varphi'_R|m_x|^2dx-\frac12\int_\mathbb{R}\varphi'_R|n|^2dx-\int_\mathbb{R}(\varphi_Rn)_x|u|^2dx.
\end{align}
Then
\begin{equation}\label{3.3}
\begin{split}
2\frac d{dt}\textup{Re}\int_\mathbb{R}\varphi_R(x-y(t))u_x\bar{u}_tdx=&-2\dot{y}\textup{Re}\int_\mathbb{R}\varphi'_R(x-y(t))u_x\bar{u}_tdx+\int_\mathbb{R}(\varphi_Rn)_x|u|^2dx\\
&-\int_\mathbb{R}\varphi'_R(x-y(t))(|u_t|^2+|u_x|^2-|u|^2)dx
\end{split}
\end{equation}
and
\begin{equation}\label{3.4}
\begin{split}
\frac1{c_0^2}\frac d{dt}\int_\mathbb{R}\varphi_R(x-y(t))nm_xdx=&-\frac1{c_0^2}\dot{y}\int_\mathbb{R}\varphi'_R(x-y(t))nm_xdx-\int_\mathbb{R}(\varphi_Rn)_x|u|^2dx\\
&-\frac1{2c_0^2}\int_\mathbb{R}\varphi'_R(x-y(t))|m_x|^2dx-\frac12\int_\mathbb{R}\varphi'_R(x-y(t))|n|^2dx.
\end{split}
\end{equation}
Adding \eqref{3.2}-\eqref{3.4} up, we get \eqref{3.5}.\hfill$\Box$

From Lemma \ref{lem3-1}, since the traditional virial identity loses positivity, we should search a modified virial quantity to overcome it. The following lemma is the key observation in the proof of the instability for the critical frequency $|\omega|=\frac1{\sqrt{2}}$.
\begin{lem}\label{lem3.2-0}
Let
\begin{align}
\tilde{I}(t)=I(t)+\textup{Re}\int_\mathbb{R} u\bar{u}_tdx-\frac1{c_0^2}\int_\mathbb{R} nmdx,
\end{align}
satisfying
 \begin{align}\label{3.8}
   \int_{\mathbb{R}}n_tdx=0,\qquad \int_{\mathbb{R}}xn_tdx=0.
 \end{align}
If $|\omega|=\frac1{\sqrt{2}}$, then
\begin{equation}\label{3.6}
\begin{split}
\tilde{I}'(t)=&-\dot{y}P(\vec{u}_0)+\dot{y}\int_\mathbb{R}[1-\varphi'_R(x-y(t))][2\textup{Re}(u_x\bar{u}_t)+\frac1{c_0^2}nm_x]dx\\
&+\int_\mathbb{R}[1-\varphi'_R(x-y(t))][|u_t|^2+|u_x|^2-|u|^2+\frac1{2c_0^2}|m_x|^2+\frac12|n|^2]dx\\
&-6E(\vec{u}_0)+8\omega Q(\vec{u}_0)+4\|u_t-i\omega u\|^2+\frac2{c_0^2}\|m_x\|^2.
\end{split}
\end{equation}
\end{lem}
{\it Proof.}  Condition \eqref{3.8} ensures the last integral term of $\tilde{I}(t)$ is well-defined in the energy space $X$. In fact, from the properties of Fourier transform, \eqref{3.8} implies
  \begin{align*}
   \int n_tdx=\int\delta(\hat{x})\hat{n_t}d\hat{x}=\hat{n_t}|_{\hat{x}=0}=0
 \end{align*}
 and
   \begin{align*}
   \int xn_tdx=i\int\delta(\hat{x})\partial_{\hat{x}}\hat{n_t}d\hat{x}=i\partial_{\hat{x}}\hat{n_t}|_{\hat{x}=0}=0.
 \end{align*}
 So the integral
 \begin{align*}
   \int nm dx=\int \hat{n}\hat{m}d\hat{x}=\int \frac{\hat{n}\hat{n_t}}{|\hat{x}|^2}d\hat{x}
 \end{align*}
 has no singularity at $\hat{x}=0$. Here we have used the Parseval identity.

By equation \eqref{1.2}, we discover that
\begin{equation}\label{3.9-0}
  \begin{split}
    \|n\|^2+\int_\mathbb{R}n|u|^2dx&=\int_\mathbb{R}n(n+|u|^2)dx\\
    &=\frac1{c_0^2}\int_\mathbb{R}n\partial_tmdx\\
    &=\frac1{c_0^2}\frac d{dt}\int_\mathbb{R} nmdx-\frac1{c_0^2}\int_\mathbb{R}\partial_x^2mmdx\\
    &=\frac1{c_0^2}\frac d{dt}\int_\mathbb{R} nmdx+\frac1{c_0^2}\|m_x\|^2
  \end{split}
\end{equation}
and
\begin{equation}
  \begin{split}
   &\|u_t\|^2-3\|u_x\|^2+(4\omega^2-3)\|u\|^2-\frac32\|n\|^2-\frac3{2c_0^2}\|m_x\|^2-3\int_\mathbb{R}n|u|^2dx\\
   =&4\|u_t-i\omega u\|^2-6E(\vec{u}_0)+8\omega Q(\vec{u}_0).
  \end{split}
\end{equation}
Thus, it follows from \eqref{3.5} and \eqref{3.2} that
\begin{equation*}
\begin{split}
\tilde{I}'(t)=&-\dot{y}P(\vec{u}_0)+\dot{y}\int_\mathbb{R}[1-\varphi'_R(x-y(t))][2\textup{Re}(u_x\bar{u}_t)+\frac1{c_0^2}nm_x]dx\\
&+\int_\mathbb{R}[1-\varphi'_R(x-y(t))][|u_t|^2+|u_x|^2-|u|^2+\frac1{2c_0^2}|m_x|^2+\frac12|n|^2]dx\\
&+\|u_t\|^2-3\|u_x\|^2-\|u\|^2-\frac32\|n\|^2+\frac1{2c_0^2}\|m_x\|^2-3\int_\mathbb{R}n|u|^2dx\\
=&-\dot{y}P(\vec{u}_0)+\dot{y}\int_\mathbb{R}[1-\varphi'_R(x-y(t))][2\textup{Re}(u_x\bar{u}_t)+\frac1{c_0^2}nm_x]dx\\
&+\int_\mathbb{R}[1-\varphi'_R(x-y(t))][|u_t|^2+|u_x|^2-|u|^2+\frac1{2c_0^2}|m_x|^2+\frac12|n|^2]dx\\
&+\|u_t\|^2-3\|u_x\|^2+(4\omega^2-3)\|u\|^2-\frac32\|n\|^2-\frac3{2c_0^2}\|m_x\|^2-3\int_\mathbb{R}n|u|^2dx\\
&+(2-4\omega^2)\|u\|^2+\frac2{c_0^2}\|m_x\|^2\\
=&-\dot{y}P(\vec{u}_0)+\dot{y}\int_\mathbb{R}[1-\varphi'_R(x-y(t))][2\textup{Re}(u_x\bar{u}_t)+\frac1{c_0^2}nm_x]dx\\
&+\int_\mathbb{R}[1-\varphi'_R(x-y(t))][|u_t|^2+|u_x|^2-|u|^2+\frac1{2c_0^2}|m_x|^2+\frac12|n|^2]dx\\
&-6E(\vec{u}_0)+8\omega Q(\vec{u}_0)+4\|u_t-i\omega u\|^2+(2-4\omega^2)\|u\|^2+\frac2{c_0^2}\|m_x\|^2.
\end{split}
\end{equation*}
In particular, if $|\omega|=\frac1{\sqrt{2}}$, we get \eqref{3.6}.\hfill$\Box$

\section{Modulation and Some Estimates}
To prove the main theorem by a contradiction argument, we suppose that the standing wave $(u_\omega,n_\omega)$ of the KGZ system is stable. The contradiction will be derived from the control of the terms in the modified virial identity. In the process of control, the use of the coercivity property would require to deal with the orthogonality conditions. The following modulation lemma shows that if the standing wave
solution is stable, then the orthogonality conditions can be obtained.
%
\begin{lem}[Modulation]
There exists $\varepsilon_0>0$ such that for any given $\varepsilon\in(0,\varepsilon_0)$, if $\vec{u}\in U_\varepsilon(\vec{\Phi}_\omega)$, then there exist $C^1$-functions
$$\theta(t):[0,t^*]\rightarrow\mathbb{R},\quad \lambda(t):[0,t^*]\rightarrow (0,\infty),\quad y(t):[0,t^*]\rightarrow\mathbb{R},$$
such that if we define
  $\vec{\xi}=(\xi,\eta,\zeta,\iota)$ satisfying
\begin{align}\label{4.1}
\vec{\xi}(t)=\vec{u}(t,\cdot-y(t))-T(\theta(t))\vec{\Phi}_{\lambda\omega}(\cdot-y(t)):=\vec{u}-\vec{R},
\end{align}
then for any given $t\in[0,t^*]$, $\vec{\xi}$ satisfies the following orthogonality conditions,
\begin{align}\label{4.5}
\langle \vec{\xi},T(\theta(t))\vec{\Upsilon}_{\lambda(t)\omega}\rangle=\langle\vec{\xi},T(\theta(t))\partial_x\vec{\Phi}_{\lambda(t)\omega}\rangle=\langle\vec{\xi},T(\theta(t))\vec{\Psi}_{\lambda(t)\omega}\rangle=0,
\end{align}
where
\begin{align*}
  T(\theta(t))\vec{\Upsilon}_{\lambda\omega}&=\Big(ie^{i\theta(t)}\phi_{\lambda\omega}(\cdot-y(t)),-\lambda\omega e^{i\theta(t)}\phi_{\lambda\omega}(\cdot-y(t)),0,0\Big)^T,\\ T(\theta(t))\partial_x\vec{\Phi}_{\lambda\omega}&=\Big(e^{i\theta(t)}\partial_x\phi_{\lambda\omega}(\cdot-y(t)),i\lambda\omega e^{i\theta(t)}\partial_x\phi_{\lambda\omega}(\cdot-y(t)),-2\phi_{\lambda\omega}\partial_x\phi_{\lambda\omega}(\cdot-y(t)),0\Big)^T,\\
  T(\theta(t))\vec{\Psi}_{\lambda\omega}&=\Big(2\lambda\omega e^{i\theta(t)}\phi_{\lambda\omega}(\cdot-y(t)),0,0,0\Big)^T.
\end{align*}
Here $\vec{R}=T(\theta(t))\vec{\Phi}_{\lambda\omega}(\cdot-y(t))$ is called the modulated soliton.

Moreover, we have estimates
\begin{align}
\|\vec{\xi}\|_X+|\lambda-1|\lesssim \varepsilon
\end{align}
and
\begin{align}\label{3.18}
  |\dot{\theta}(t)-\lambda(t)\omega|+|\dot{y}|+|\dot{\lambda}|=O(\|\vec{\xi}\|_X).
\end{align}
\end{lem}
{\it Proof.} We first show the existence of $\theta,\ \lambda$ and $y$ by the Implicit Function Theorem. Denoting
\begin{align*}
  F_1(\theta,y,\lambda;\vec{u})&=\langle \vec{\xi},T(\theta)\vec{\Upsilon}_{\lambda\omega}\rangle,\\
   F_2(\theta,y,\lambda;\vec{u})&=\langle\vec{\xi},T(\theta)\partial_x\vec{\Phi}_{\lambda\omega}\rangle,\\
   F_3(\theta,y,\lambda;\vec{u})&=\langle\vec{\xi},T(\theta)\vec{\Psi}_{\lambda\omega}\rangle,
\end{align*}
we have
\begin{align*}
  F_1(\theta,y,\lambda;\vec{u})=&\Big\langle u-e^{i\theta}\phi_{\lambda\omega}(\cdot-y),ie^{i\theta}\phi_{\lambda\omega}(\cdot-y)\Big\rangle+\Big\langle v-i\lambda\omega e^{i\theta}\phi_{\lambda\omega}(\cdot-y),-\lambda\omega e^{i\theta}\phi_{\lambda\omega}(\cdot-y)\Big\rangle,\\
   F_2(\theta,y,\lambda;\vec{u})=&\Big\langle u-e^{i\theta}\phi_{\lambda\omega}(\cdot-y),e^{i\theta}\partial_x\phi_{\lambda\omega}(\cdot-y)\Big\rangle+\Big\langle v-i\lambda\omega e^{i\theta}\phi_{\lambda\omega}(\cdot-y),i\lambda\omega e^{i\theta}\partial_x\phi_{\lambda\omega}(\cdot-y)\Big\rangle\\
   &+\Big\langle n+|\phi_{\lambda\omega}(\cdot-y)|^2,-2\phi_{\lambda\omega}\partial_x\phi_{\lambda\omega}(\cdot-y)\Big\rangle,\\
   F_3(\theta,y,\lambda;\vec{u})=&\Big\langle u-e^{i\theta}\phi_{\lambda\omega}(\cdot-y),2\lambda\omega e^{i\theta}\phi_{\lambda\omega}(\cdot-y)\Big\rangle,
\end{align*}
then
\begin{equation*}
  \frac{\partial(F_1,F_2,F_3)}{\partial(\theta,y,\lambda)}\Big|_{(0,0,1;\vec{\Phi}_\omega)}=\left(
                                                             \begin{array}{ccc}
                                                               a_{11} & 0 & 0 \\
                                                               0 & a_{22} & a_{23} \\
                                                               0 & 0 & a_{33} \\
                                                             \end{array}
                                                           \right),
\end{equation*}
where
\begin{align*}
 &a_{11}=\partial_\theta F_1\big|_{(0,0,1;\vec{\Phi}_\omega)}=-(1+\omega^2)\|\phi_\omega\|^2,\\
 & a_{22}=\partial_y F_2\big|_{(0,0,1;\vec{\Phi}_\omega)}=(1+\omega^2)\|\partial_x\phi_\omega\|^2+\|\partial_x|\phi_\omega|^2\|^2,\\
 &  a_{23}=\partial_\lambda F_2\big|_{(0,0,1;\vec{\Phi}_\omega)}=\frac{\omega^2}{1-\omega^2}[(1+\omega^2)\|x^\frac12\partial_x\phi_\omega\|^2+\|x^\frac12\partial_x|\phi_\omega|^2\|^2],\\
 &a_{33}= \partial_\lambda F_3\big|_{(0,0,1;\vec{\Phi}_\omega)}=\frac{\omega^3}{1-\omega^2}\|\phi_\omega\|^2.
\end{align*}
This means that the Jacobian matrix of the derivative of $(\theta,y,\lambda)\rightarrow (F_1,F_2,F_3)$ is nondegenerate at $(\theta,y,\lambda)=(0,0,1)$. Hence, by the Implicit Function Theorem, we obtain the existence of $(\theta,y,\lambda)$ in a neighborhood of $(0,0,1)$ satisfying
$$(F_1,F_2,F_3)(\theta(t),y(t),\lambda(t);\vec{u}(t))=0.$$
To verify the modulation parameters $(\theta,y,\lambda)$ are $C^1$, we can use the equation of $\vec{\xi}$ and regularization arguments, similar to the method in
\cite{martel}.

We now want to verify \eqref{3.18}. 
From \eqref{4.1} and \eqref{3.16}, we have
\begin{equation}
  \partial_t\vec{\xi}+i\dot{\theta}\vec{R}-\dot{y}\partial_x \vec{R}+\dot{\lambda}\partial_\lambda\vec{R}=JE'(\vec{R}+\vec{\xi}).
\end{equation}
Moreover, from Taylor's expansion, we have
\begin{equation}
 E'(\vec{R}+\vec{\xi})=E'(\vec{R})+E''(\vec{R})\vec{\xi}+O(\|\vec{\xi}\|_X^2).
\end{equation}
Noting the fact that $\phi_{\lambda\omega}$ is a solution of
\begin{equation}
  -\partial_x^2\phi_{\lambda\omega}+(1-\lambda^2\omega^2)\phi_{\lambda\omega}-\phi^3_{\lambda\omega}=0,
\end{equation}
we have
\begin{align}\label{3.20}
  \partial_t\vec{\xi}+i(\dot{\theta}-\lambda\omega)A\vec{R}-\dot{y}\partial_x\vec{R}+\dot{\lambda}\partial_\lambda\vec{R}=JE''(\vec{R})\vec{\xi}+O(\|\vec{\xi}\|_X^2),
\end{align}
where
$$A=\begin{pmatrix}
                            1  & 0 & 0 & 0 \\
                            0 & 1 & 0 & 0 \\
                            0 & 0& 0 & 0 \\
                            0 & 0 & 0 & 0 \\
                          \end{pmatrix}.$$
We differentiate $\langle\vec{\xi},T(\theta)\vec{\Psi}_{\lambda(t)\omega}\rangle=0$ in \eqref{4.5} with respect to time $t$ to get
\begin{align}
\langle\partial_t\vec{\xi},T(\theta)\vec{\Psi}_{\lambda(t)\omega}\rangle=-\langle\vec{\xi},\partial_t\big(T(\theta)\vec{\Psi}_{\lambda(t)\omega}\big)\rangle.
\end{align}
Denoting
$$\textup{Mod}(t)=(\dot{\theta}(t)-\lambda(t)\omega,\dot{y}(t),\dot{\lambda})^T,$$
we have
\begin{align}\label{3.24}
 \langle\vec{\xi},\partial_t(T(\theta)\vec{\Psi}_{\lambda(t)\omega})\rangle=O\big((1+|\textup{Mod}(t)|)\|\vec{\xi}\|_X\big).
\end{align}
Because of the orthogonality condition and taking the inner product of \eqref{3.20} with $T(\theta)\vec{\Psi}_{\lambda(t)\omega}$, we get
\begin{equation}\label{3.23}
  \begin{split}
   &\langle\partial_t\vec{\xi},T(\theta)\vec{\Psi}_{\lambda(t)\omega}\rangle\\
   =&-(\dot{\theta}-\lambda\omega)\langle iA\vec{R},T(\theta)\vec{\Psi}_{\lambda(t)\omega}\rangle+\dot{y}\langle\partial_x\vec{R},T(\theta)\vec{\Psi}_{\lambda(t)\omega}\rangle-\dot{\lambda}\langle\partial_\lambda\vec{R},T(\theta)\vec{\Psi}_{\lambda(t)\omega}\rangle+O(\|\vec{\xi}\|_X^2)\\
   =&-\dot{\lambda}\langle\partial_\lambda\phi_{\lambda\omega},2\lambda\omega\phi_{\lambda\omega}\rangle+O(\|\vec{\xi}\|_X^2).
  \end{split}
  \end{equation}
\eqref{3.24} and \eqref{3.23} imply that
\begin{align}\label{3.25}
\lambda\omega\partial_\lambda\|\phi_{\lambda\omega}\|^2\dot{\lambda}=O\big((1+|\textup{Mod}(t)|)\|\vec{\xi}\|_X\big)+O(\|\vec{\xi}\|_X^2).
\end{align}
Taking the inner product of \eqref{3.20} with $T(\theta)\vec{\Upsilon}_{\lambda\omega}$ and $T(\theta)\partial_x\vec{\Phi}_{\lambda\omega}$, respectively, by similar arguments we get
\begin{equation}\label{3.26}
  \begin{split}
    (1+\lambda^2\omega^2)\|\phi_{\lambda\omega}\|^2(\dot{\theta}-\lambda\omega)=O\big((1+|\textup{Mod}(t)|)\|\vec{\xi}\|_X\big)+O(\|\vec{\xi}\|_X^2),
  \end{split}
\end{equation}
\begin{equation}\label{3.27}
  \begin{split}
    [(1+\lambda^2\omega^2)\|\partial_x\phi_{\lambda\omega}\|^2+\|\partial_x|\phi_{\lambda\omega}|^2\|^2]\dot{y}+C_\star\dot{\lambda} =O\big((1+|\textup{Mod}(t)|)\|\vec{\xi}\|_X\big)+O(\|\vec{\xi}\|_X^2).
  \end{split}
\end{equation}
Indeed, we have
\begin{equation}\label{3.28}
  \begin{split}
   M \textup{Mod}(t)=O\big((1+|\textup{Mod}(t)|)\|\vec{\xi}\|_X\big)+O(\|\vec{\xi}\|_X^2),
  \end{split}
\end{equation}
where
\begin{align*}
M=\begin{pmatrix}
    (1+\lambda^2\omega^2)\|\phi_{\lambda\omega}\|^2 & 0 & 0 \\
    0 & (1+\lambda^2\omega^2)\|\partial_x\phi_{\lambda\omega}\|^2+\|\partial_x|\phi_{\lambda\omega}|^2\|^2 & C_\star \\
    0 & 0 & \lambda\omega\partial_\lambda\|\phi_{\lambda\omega}\|^2 \\
  \end{pmatrix}
\end{align*}
is an invertible matrix. Therefore,
\begin{equation}\label{3.29}
  \begin{split}
   |\textup{Mod}(t)|\leq C\|\vec{\xi}\|_X+O(\|\vec{\xi}\|_X^2),
  \end{split}
\end{equation}
which concludes the proof of \eqref{3.18}.\hfill$\Box$

In what follows, we want to scale the terms on the righthand side of the modified virial identity \eqref{3.6}.
For this purpose, we first choose the initial data as
\begin{align}\label{3.9}
  \vec{u}_0=(1+a)\vec{\Phi}_\omega
\end{align}
and get the following properties.
\begin{lem}\label{lem3.2}
Let $\vec{u}_0$ be defined by \eqref{3.9}. For the KGZ system, we have
  \begin{align}
    P(\vec{u}_0)&=0,\label{3.10}\\
    Q(\vec{u}_0)-Q(\vec{\Phi}_\omega)&=2a\omega\|\phi_\omega\|^2+O(a^2),\label{3.11}\\
    Q(T(\theta)\vec{\Phi}_\omega)-Q(\vec{\Phi}_\omega)&=0\label{3.11-0}
  \end{align}
  and
  \begin{align}\label{3.12}
    -6E(\vec{u}_0)+8\omega Q((\vec{u}_0))=(4a\omega^2 +4\omega^2-2)\|\phi_\omega\|^2+O(a^2),
  \end{align}
where $O(a^2)$ denotes the equivalent infinitesimal of $a^2$.
\end{lem}
\textit{Proof.} To be specific,
$$\vec{u}_0=\big((1+a)\phi_\omega,i(1+a)\omega\phi_\omega,-(1+a)\phi_\omega^2,0\big),$$
then
$$P(\vec{u}_0)=2\textup{Re}\int i(1+a)^2\omega\phi_\omega\partial_x\phi_\omega dx+0=0,$$
which gives \eqref{3.10}.

 For getting \eqref{3.11}, we have
\begin{align*}
  Q(\vec{u}_0)- Q(\vec{\Phi}_\omega)=&\langle Q'(\vec{\Phi}_\omega),\vec{u}_0-\vec{\Phi}_\omega\rangle+O(\|\vec{u}_0-\vec{\Phi}_\omega\|_{X}^2) \\
=&\langle \omega\phi_\omega,a\phi_\omega\rangle+\langle i\phi_\omega,ia\omega\phi_\omega\rangle+O(a^2)\\
=&2a\omega\|\phi_\omega\|^2+O(a^2).
\end{align*}
Property \eqref{3.11-0} is a direct consequence from the definition of the charge. As to \eqref{3.12}, we write
\begin{align*}
  -6E(\vec{u}_0)+8\omega Q((\vec{u}_0))=&-6\big(E(\vec{u}_0)-E(\vec{\Phi}_\omega)\big)+8\omega\big(Q(\vec{u}_0)-Q(\vec{\Phi}_\omega)\big)\\
  & -6E(\vec{\Phi}_\omega)+8\omega Q((\vec{\Phi}_\omega))\\
=&-6\big(S_\omega(\vec{u}_0)-S_\omega(\vec{\Phi}_\omega)\big)+2\omega\big(Q(\vec{u}_0)-Q(\vec{\Phi}_\omega)\big)\\
  & -6E(\vec{\Phi}_\omega)+8\omega Q((\vec{\Phi}_\omega)).
\end{align*}
According to Taylor's expansion and \eqref{2.4-0}, we have
$$S_\omega(\vec{u}_0)-S_\omega(\vec{\Phi}_\omega)=O\big(\|\vec{u}_0-\vec{\Phi}_\omega\|_X^2\big)=O(a^2).$$
Thus, combined with \eqref{3.7} in Lemma \ref{lem3.1}, we get \eqref{3.12}.\hfill$\Box$

Moreover, we have the following properties.
\begin{lem}\label{lem4.4}
Let $\vec{u}_0$ be defined by \eqref{3.9}, $\lambda\in\mathbb{R}^+$ and $\lambda\lesssim 1$. For the KGZ system, we have
  \begin{align}
    S_{\lambda\omega}(\vec{u}_0)-S_{\lambda\omega}(\vec{\Phi}_\omega)&=-2a(\lambda-1)\omega^2\|\phi_\omega\|^2+O(a^2),\label{3.12-0}\\
    S_{\lambda\omega}(T(\theta)\vec{\Phi}_\omega)-S_{\lambda\omega}(\vec{\Phi}_\omega)&=0\label{4.22}
  \end{align}
  and
 \begin{equation}\label{lem5.1}
   S_{\lambda\omega}(T(\theta)\vec{\Phi}_{\lambda\omega})-S_{\lambda\omega}(\vec{\Phi}_{\omega})=-\frac{(\lambda-1)^2(1-2\omega^2)(2-\omega)\omega^2}{2(1-\omega^2)}\|\phi_\omega\|^2+o(|\lambda-1|^2).
  \end{equation}
\end{lem}
\textit{Proof.} From Lemma \ref{lem3.2}, we have
\begin{align*}
  S_{\lambda\omega}(\vec{u}_0)-S_{\lambda\omega}(\vec{\Phi}_\omega)&=E(\vec{u}_0)-E(\vec{\Phi}_{\omega})-\lambda\omega Q(\vec{u}_0)+\lambda\omega Q(\vec{\Phi}_{\omega})\\
  &=S_\omega(\vec{u}_0)-S_\omega(\vec{\Phi}_{\omega})-(\lambda-1)\omega [Q(\vec{u}_0)- Q(\vec{\Phi}_{\omega})]\\
  &=O(a^2)-2a(\lambda-1)\omega^2\|\phi_\omega\|^2,
\end{align*}
which gives \eqref{3.12-0}.

Since
  \begin{align*}
    S_{\lambda\omega}(T(\theta)\vec{\Phi}_\omega)-S_{\lambda\omega}(\vec{\Phi}_\omega)=E(T(\theta)\vec{\Phi}_\omega)-E(\vec{\Phi}_{\omega})-\lambda\omega[ Q(T(\theta)\vec{\Phi}_\omega)- Q(\vec{\Phi}_{\omega})],
  \end{align*}
 \eqref{4.22} can be checked by the definition of the energy and \eqref{3.11-0}.

From the definition of $S_\omega$, Taylor's expansion, \eqref{2.4-0}  and \eqref{4.22}, we have
\begin{equation}
\begin{split}
  &S_{\lambda\omega}(T(\theta)\vec{\Phi}_{\lambda\omega})-S_{\lambda\omega}(\vec{\Phi}_{\omega})\\
  =&S_{\lambda\omega}(T(\theta)\vec{\Phi}_{\lambda\omega})-S_{\lambda\omega}(T(\theta)\vec{\Phi}_{\omega})+S_{\lambda\omega}(T(\theta)\vec{\Phi}_{\omega})-S_{\lambda\omega}(\vec{\Phi}_{\omega})\\
  =&S_{\omega}(T(\theta)\vec{\Phi}_{\lambda\omega})-S_{\omega}(T(\theta)\vec{\Phi}_{\omega})+(\lambda-1)\omega[Q(T(\theta)\vec{\Phi}_{\lambda\omega})-Q(T(\theta)\vec{\Phi}_{\omega})]\\
  =&\frac12\langle S_\omega''(\vec{\Phi}_{\omega})\big(T(\theta)\vec{\Phi}_{\lambda\omega}-T(\theta)\vec{\Phi}_{\omega}\big),\big(T(\theta)\vec{\Phi}_{\lambda\omega}-T(\theta)\vec{\Phi}_{\omega}\big)\rangle+o(|\lambda-1|^2)\\
  &+(\lambda-1)\omega[Q(T(\theta)\vec{\Phi}_{\lambda\omega})-Q(T(\theta)\vec{\Phi}_{\omega})].
\end{split}
\end{equation}
Since
\begin{align*}
  \langle &S_\omega''(\vec{\Phi}_{\omega})\big(T(\theta)\vec{\Phi}_{\lambda\omega}-T(\theta)\vec{\Phi}_{\omega}\big),\big(T(\theta)\vec{\Phi}_{\lambda\omega}-T(\theta)\vec{\Phi}_{\omega}\big)\rangle\\
  =&(\lambda-1)^2\omega^2\langle S_\omega''(T(\theta)\vec{\Phi}_{\omega})\partial_\omega(T(\theta)\vec{\Phi}_\omega),\partial_\omega(T(\theta)\vec{\Phi}_\omega)\rangle+o(|\lambda-1|^2)\\
  =&-(\lambda-1)^2\omega^2\frac d{d\lambda}Q(T(\theta)\vec{\Phi}_{\lambda\omega})\big|_{\lambda=1}+o(|\lambda-1|^2)
\end{align*}
and
$$\frac d{d\lambda}Q(T(\theta)\vec{\Phi}_{\lambda\omega})\big|_{\lambda=1}=\omega \frac d{d\omega}Q(\vec{\Phi}_{\omega})=-\frac{(1-2\omega^2)\omega}{1-\omega^2}\|\phi_\omega\|^2,$$
we have
\begin{align*}
  \langle &S_\omega''(\vec{\Phi}_{\omega})\big(T(\theta)\vec{\Phi}_{\lambda\omega}-T(\theta)\vec{\Phi}_{\omega}\big),\big(T(\theta)\vec{\Phi}_{\lambda\omega}-T(\theta)\vec{\Phi}_{\omega}\big)\rangle\\
  =&\frac{(\lambda-1)^2(1-2\omega^2)\omega^3}{1-\omega^2}\|\phi_\omega\|^2+o(|\lambda-1|^2).
\end{align*}
Due to
\begin{align*}
 Q(T(\theta)\vec{\Phi}_{\lambda\omega})-Q(T(\theta)\vec{\Phi}_{\omega})=&(\lambda-1)\omega\langle Q'(T(\theta)\vec{\Phi}_{\omega}),\partial_\omega (T(\theta)\vec{\Phi}_{\omega})\rangle +o(|\lambda-1|)\\
 =&(\lambda-1)\omega\frac d{d\omega}Q(\vec{\Phi}_{\omega})+o(|\lambda-1|)\\
 =&-\frac{(\lambda-1)(1-2\omega^2)\omega}{1-\omega^2}\|\phi_\omega\|^2+o(|\lambda-1|),
\end{align*}
we get \eqref{lem5.1} immediately.\hfill$\Box$

Now we are ready to control $\|u_t-i\omega u\|$.
\begin{lem}\label{lem4.2}
Suppose that $\vec{\xi}=(\xi,\eta,\zeta,0)$ is defined by \eqref{4.1}. If $|\omega|=\frac1{\sqrt{2}}$, then
  \begin{align}\label{4.6}
   \|u_t-i\omega u\|^2=(\lambda-1)^2\omega^2\|\phi_{\omega}\|^2+\|\eta-i\omega\xi\|^2+O\big(a|\lambda-1|+|\lambda-1|^3+|\lambda-1|\|\vec{\xi}\|^2_X\big).
  \end{align}
\end{lem}
\textit{Proof.} By \eqref{4.1}, we have
\begin{align*}
\|u_t-i\omega u\|^2=&\|ie^{i\theta}\lambda\omega\phi_{\lambda\omega}+\eta-i\omega e^{i\theta}\phi_{\lambda\omega}-i\omega\xi\|^2\\
=&\|ie^{i\theta}(\lambda-1)\omega\phi_{\lambda\omega}+\eta-i\omega\xi\|^2\\
=&(\lambda-1)^2\omega^2\|\phi_{\lambda\omega}\|^2+\|\eta-i\omega\xi\|^2+2(\lambda-1)\omega\langle i e^{i\theta}\phi_{\lambda\omega},\eta-i\omega\xi\rangle.
\end{align*}
Since the third orthogonality condition $\langle\vec{\xi},2\lambda\omega e^{i\theta}\phi_{\lambda\omega}\rangle=0$ in \eqref{4.5} and
$$\|\phi_{\lambda\omega}\|^2=\|\phi_{\omega}\|^2+O(|\lambda-1|),$$
we get
\begin{align*}
\|u_t-i\omega u\|^2=(\lambda-1)^2\omega^2\|\phi_{\omega}\|^2+\|\eta-i\omega\xi\|^2+2(\lambda-1)\omega\langle i e^{i\theta}\phi_{\lambda\omega},\eta\rangle+O(|\lambda-1|^3).
\end{align*}
To estimate $\langle i e^{i\theta}\phi_{\lambda\omega},\eta\rangle$, from \eqref{3.11}, \eqref{3.11-0} and the charge conservation law we have
\begin{align*}
  &Q(\vec{u}_0)-Q(\vec{\Phi}_\omega)+[Q(\vec{\Phi}_\omega)-Q(T(\theta)\vec{\Phi}_{\omega})]+[Q(T(\theta)\vec{\Phi}_\omega)-Q(T(\theta)\vec{\Phi}_{\lambda\omega})]\\
  =&Q(\vec{u})-Q(T(\theta)\vec{\Phi}_{\lambda\omega})\\
  =&-\langle \xi, \lambda\omega e^{i\theta}\phi_{\lambda\omega}\rangle+\langle \eta, i e^{i\theta}\phi_{\lambda\omega}\rangle+O(\|\vec{\xi}\|^2_X).
\end{align*}
Recalling
 $\frac d{dw}Q(\vec{\Phi}_\omega)=0$ when $|\omega|=\frac1{\sqrt{2}}$ in \eqref{3.6-0}, we have
\begin{align*}
  \langle\eta,ie^{i\theta}\phi_{\lambda\omega}\rangle=&Q(\vec{u}_0)-Q(\vec{\Phi}_\omega)+[Q(T(\theta)\vec{\Phi}_\omega)-Q(T(\theta)\vec{\Phi}_{\lambda\omega})]+O(\|\vec{\xi}\|^2_X)\\
  =&O(a+|\lambda-1|^2+\|\vec{\xi}\|^2_X).
\end{align*}
Thus, \eqref{4.6} is obtained.
\hfill$\Box$


We then scale $\|\vec{\xi}\|_X$.
\begin{lem}\label{lem5.2}
  Let $\vec{\xi}$ be defined by \eqref{4.1}. If $|\omega|=\frac1{\sqrt{2}}$, then
  \begin{align}\label{5.3}
   \|\vec{\xi}\|_X^2=O(a^2+a|\lambda-1|)+o(|\lambda-1|^2).
  \end{align}
\end{lem}
\textit{Proof.} By \eqref{2.4-0} and Lemma \ref{lem2.3}, combined with the energy and charge conservation laws, we have
\begin{align*}
  \|\vec{\xi}\|_X^2\lesssim& \langle S''_{\lambda\omega}(\vec{\Phi}_{\lambda\omega})\vec{\xi},\vec{\xi}\rangle\\
  =&S_{\lambda\omega}(\vec{u})-S_{\lambda\omega}(T(\theta)\vec{\Phi}_{\lambda\omega})+o(\|\vec{\xi}\|_X^2)\\
  =&[S_{\lambda\omega}(\vec{u}_0)-S_{\lambda\omega}(\vec{\Phi}_{\omega})]-[S_{\lambda\omega}(T(\theta)\vec{\Phi}_{\lambda\omega})-S_{\lambda\omega}(\vec{\Phi}_{\omega})]+o(\|\vec{\xi}\|_X^2).
\end{align*}
When $|\omega|=\frac1{\sqrt{2}}$, from \eqref{lem5.1} we have
  \begin{equation}
   S_{\lambda\omega}(T(\theta)\vec{\Phi}_{\lambda\omega})-S_{\lambda\omega}(\vec{\Phi}_{\omega})=o(|\lambda-1|^2),
  \end{equation}
 then, thanks to \eqref{3.12-0}, we get \eqref{5.3} immediately.\hfill$\Box$

\section{Proof of Theorem \ref{thm}}
 Suppose that $\vec{u}\in U_\varepsilon(\vec{\Phi}_\omega)$ defined by \eqref{4.1}, we have $|\lambda-1|\lesssim \varepsilon<<1$. We will prove Theorem \ref{thm} by a contradiction argument. The initial data satisfies \eqref{3.9} in which $a>0$ will be determined later. Assume that $\vec{\xi}(t,x)=(\xi,\eta,\zeta,\iota)^T$ satisfies
 \begin{align}
   \int_{\mathbb{R}}\iota dx=0,\qquad \int_{\mathbb{R}}x\cdot\iota dx=0,
 \end{align}
  by the definition of $\tilde{I}(t)$ and the proof of Lemma \ref{lem3.2-0}, we have the time uniform boundedness
\begin{align*}
  \tilde{I}(t)\lesssim R(\|\vec{\Phi}_\omega\|^2_X+1).
\end{align*}
%
However, we will claim that
\begin{align}\label{5.2}
\lim_{t\rightarrow\infty} \tilde{I}(t)\rightarrow\infty,
\end{align}
which gives a contradiction. Since the solution to \eqref{1.5} is exponentially decaying at infinity, we have
\begin{align*}
&\int_{|x-y(t)|\geq R}[|u_t|^2+|u_x|^2-|u|^2+\frac1{2c_0^2}|m_x|^2+\frac12|n|^2+2\textup{Re}(u_x\bar{u}_t)+\frac1{c_0^2}nm_x]dx\\
\lesssim&\int_{|x|\geq R}(|\phi_{\lambda\omega}|^2+|\partial_x\phi_{\lambda\omega}|^2+|\xi|^2+|\xi_x|^2+|\eta|^2+|(-\Delta)^{-1}\iota_x|^2+|\zeta|^2)dx\\
=&O\big(\|\vec{\xi}\|^2_X+\frac1R\big).
\end{align*}
When $|\dot{y}|\lesssim 1$, from Lemma \ref{lem3.2} we have
\begin{align*}
  \tilde{I}'(t)\gtrsim &(4a\omega^2 +4\omega^2-2)\|\phi_\omega\|^2+O(a^2)+4(\lambda-1)^2\omega^2\|\phi_{\lambda\omega}\|^2+4\|\eta-i\omega\xi\|^2\\
  &+O\big(a|\lambda-1|+|\lambda-1|^3+|\lambda-1|\|\vec{\xi}\|^2_X\big)+O\big(\|\vec{\xi}\|^2_X+\frac1R\big)\\
  \gtrsim &(4a\omega^2 +4\omega^2-2)\|\phi_\omega\|^2+O\big(a^2 +a|\lambda-1|+|\lambda-1|^3+\|\vec{\xi}\|^2_X+\frac1R\big).
\end{align*}
Let $R$ be sufficiently large and $\frac1R=O(a^2)$. According to Lemma \ref{lem5.2}, we have
\begin{align*}
  \tilde{I}'(t)\gtrsim &(4a\omega^2 +4\omega^2-2)\|\phi_\omega\|^2+O\big(a^2 +a|\lambda-1|\big)+o(|\lambda-1|^2).
\end{align*}
Since $|\omega|=\frac1{\sqrt{2}}$, $4a\omega^2$ is absolutely positive. Choosing $\varepsilon$ and $a_0$ small enough, for any given $a\in(0,a_0)$, we have
$$\tilde{I}'(t)\gtrsim 2a\omega^2 \|\phi_\omega\|^2>0,$$
which means \eqref{5.2}. Thus, we have completed the instability of the standing wave at frequency $|\omega|=\frac1{\sqrt{2}}$ for the KGZ system.


\section*{Acknowledgements}
The author would like to thank Professor Yifei Wu for his encouragements and Professor Yi Zhou for his suggestion on condition \eqref{3.8}. The author was supported by China Postdoctoral Science Foundation (2017M611528).


\begin{thebibliography}{99}

\bibitem{ambrosetti}
Ambrosetti A., Malchiodi A., Nonlinear Anaylsis and Semilinear Elliptic Problems. Cambridge: Cambridge University Press, 2007.

\bibitem{b}
Bellazzini J., Ghimenti M., Le Coz S., Multi-solitary waves for the nonlinear Klein-Gordon equation. Comm. Partial Differ Equ., 39(2014), pp. 1479-1522.
\bibitem{comech}
Comech A., Pelinovsky D., Purely nonlinear istability of standing waves with minimal energy. Comm. Pure Appl. Math., 56(2003), pp. 1565-1607.

\bibitem{c}
Le Coz S., Wu Y. F., Stability of multisolitons for the derivative nonlinear Schr\"{o}dinger equation, International Mathematics Research Notices, 2017, doi: 10.1093/imrn/rnx013.

\bibitem{dendy}
Dendy R. O., Plasma Dynamics, Oxford University Press, Oxford, UK, 1990.



\bibitem{gan}
Gan Z. H., Zhang, J., Instability of standing waves for Klein-Gordon-Zakharov equations with different propagation speeds in three space dimensions, J. Math. Anal. Appl., 307 (2005), pp. 219-231.

\bibitem{gan09}
Gan Z. H., Guo B. L., Zhang J., Instability of standing wave, global existence and blowup for the Klein-Gordon-Zakharov system with different-degree nonlinearities. J. Differential Equations, 246 (2009), 4097-4128.
\bibitem{gidas}
Gidas B., Ni W.M., Nirenberg L., Symmetry and related properties via the maximum principle. Comm. Math. Phys., 68 (1979), pp. 209-243.
\bibitem{glangetas}
Glangetas L., Merle F., Concentration properties of blow-up solutions and instability results for Zakharov equation in dimension two. Part II. Commun. Math. Phys. 160 (1994), pp. 349-389.

\bibitem{grillakis1}
Grillakis M., Shatah J., Strauss W., Stability theory of solitary waves in the presence of symmetry I. J. Funct. Anal., 74 (1987), pp. 160-197.

\bibitem{grillakis2}
Grillakis M., Shatah J., Strauss W., Stability theory of solitary waves in the Presence of Symmetry II. J. Funct. Anal., 94 (1990), pp. 308-348.
\bibitem{guo95}
Guo B., Yuan G. W., Global smooth solution for the Klein-Gordon-Zakharov equations, J. Math. Phys., 36 (1995), 4119-4124.

\bibitem{guozihua}
Guo Z. H., Nakanishi K., Wang S. X., Global dynamics below the ground state energy for the Klein-Gordon-Zakharov system in the 3D radial case. Comm. Partial Differential Equations, 39 (2014), pp. 1158-1184.
\bibitem{hakkaev}
Hakkaev S., Stanislavova M., Stefanov A., Orbital stability for periodic standing waves of the Klein-Gordon-Zakharov system and the beam equation. Z. Angew. Math. Phys., 64 (2013), pp. 265-282.

\bibitem{kato}
Kato I., Kinoshita S., Well-posedness for the Cauchy problem of the Klein-Gordon-Zakharov system in five and more dimensions. arXiv: 1612.04240v1.
\bibitem{kinoshita}
Kinoshita S., Well-posedness for the Cauchy problem of the Klein-Gordon-Zakharov system in 2D. arXiv: 1612.04241v1.
\bibitem{kwong}
Kwong M.K., Uniqueness of positive solutions of $\Delta u-u+u^p=0$ in $\mathbb{R}^n$. Arch. Rational Mech. Anal., 105 (1989), pp. 243-266.
\bibitem{maeda}
Maeda M., Stability of bound states of Hamiltonian PDEs in the degenerate cases. J. Funct. Anal., 263(2012), pp. 511-528.

\bibitem{martel}
Martel Y., Merle F., Instability of solitons for the critical generalized Korteweg-de Vries equation. Geom. Funct. Anal. 11 (2001), pp. 74-123.

\bibitem{mm}
Martel Y., Merle F., Tsai T. P., Stability in $H^1$ of the sum of $K$ solitary waves for some nonlinear Schr\"{o}dinger equations. Duke Mathematical Journal, 133(2006), pp. 405-466.

\bibitem{m}
Merle F. Blow-up results of virial type for Zakharov equations, Comm. Math. Phys., 175 (1996), pp. 433-455.
\bibitem{merle}
Merle F., Existence of blow-up solutions in the energy space for the critical generalized KdV
equation. J. Amer. Math. Soc., 14 (2001), pp. 555¨C578.
\bibitem{raphael}
Merle, F., and Rapha\"{e}l, P., On universality of blow-up profile for L2 critical nonlinear
Schr\"{o}dinger equation. Invent. Math., 156 (2004), pp. 565¨C672.
\bibitem{ohta}
Ohta M., Instability of bound states for abstract nonlinear Schrodinger equations. J. Funct. Anal., 261(2011), pp. 90-110.
\bibitem{ot}
Ohta M., Todorova G. Strong instability of standing wavew for the nonlinear Klein-Gordon equation and the Klein-Gordon-Zakharov system, SIAM J. Math. Anal., 38 (2007), pp. 1912-1931.

\bibitem{ott}
Ozawa T., Tsutaya K., Tsutsumi Y. Well-posedness in energy space for the Cauchy problem of the Klein-Gordon-Zakharov equations with different propagation speeds in three space dimensions, Math. Ann., 313 (1999), pp. 127-140.



\bibitem{strauss}
Strauss W., Existence of solitary waves in higher dimensions, Comm. Math. Phys., 55 (1977), pp. 149-162.
\bibitem{stuart}
Stuart D. M. A., modulational approach to stability of non-topological solitons in semilinear wave equations, J. Math. Pures Appl., 80 (2001), pp. 51083.
\bibitem{weinstein}
Weinstein M., Modulational stability of ground states of nonlinear Schr\"{o}dinger equations,
SIAM J. Math. Anal. 16 (1985), pp. 472¨C491.
\bibitem{w}
Wu Y. F. Instability of the standing waves for the nonlinear Klein-Gordon equations in one dimension, ArXiv: 1705.04216v2.

\bibitem{zakharov}
Zakharov V. E., Collapse of Langmuir waves, Soviet Phys. JETP, 35 (1972), pp. 908-914.

\bibitem{local}
Zhou Y., Guo B., Periodic boundary problem and initial value problem for the generalized Korteweg-de Vries systems of higher order (in Chinese), Acta Math. Sinica, 27 (1984), pp. 154-176.

\end{thebibliography}
\end{document}